\documentclass[11pt]{article}
\usepackage{amsfonts, amssymb, amsmath, amsthm, xcolor, graphicx}
\begin{document}

\newtheorem{thm}{Theorem}[section]
\newtheorem{dfn}[thm]{Definition}
\newtheorem{lm}[thm]{Lemma}
\newtheorem{prop}[thm]{Proposition}
\newtheorem{cor}[thm]{Corollary}
\newtheorem{con}[thm]{Conjecture}
\newcommand{\Fix}{\mathop{\mathrm{Fix}}}
\newcommand{\base}{\mathbf{b}}
\newcommand{\fix}{\mathop{\mathrm{fix}}}
\newcommand{\Sym}{\mathop{\mathrm{Sym}}}
\newcommand{\Alt}{\mathop{\mathrm{Alt}}}
\newcommand{\cov}{\mathop{\mathrm{cr}}}
\newcommand{\ov}{\overline}
\newcommand{\F}{\mathcal{F}}
\newcommand{\V}{\mathcal{V}}
\newcommand{\A}{\mathcal{A}}
\newcommand{\B}{\mathcal{B}}
\newcommand{\C}{\mathcal{C}}
\newcommand{\D}{\mathcal{D}}
\newcommand{\E}{\mathcal{E}}
\newcommand{\G}{\mathcal{G}}
\newcommand{\Hc}{\mathcal{H}}
\newcommand{\R}{\mathcal{R}}
\newcommand{\T}{\mathcal{T}}
\newcommand{\U}{\mathcal{U}}
\newcommand{\Se}{\mathcal{S}}
\newcommand{\Po}{\mathcal{P}}
\newcommand{\prob}{\mathbb{P}}
\newcommand{\sh}{\triangleleft}
\newcommand{\QED}{\hfill \blacksquare}

\renewcommand{\baselinestretch}{1.1}

\title{\bf Group Marriage Problem}
\author{
Cheng Yeaw Ku
\thanks{ Department of Mathematics, National University of Singapore, Singapore 117543. E-mail: matkcy@nus.edu.sg} \and K.B. Wong \thanks{
Institute of Mathematical Sciences, University of Malaya, 50603 Kuala Lumpur, Malaysia. E-mail:
kbwong@um.edu.my.} } \maketitle

\begin{abstract}\noindent
Let $G$ be a permutation group acting on $[n]=\{1, \ldots, n\}$ and $\mathcal{V}=\{V_{i}: i=1, \ldots, n\}$ be a system of $n$ subsets of $[n]$. When is there an element $g \in G$ so that $g(i) \in V_{i}$ for each $i \in [n]$? If such $g$ exists, we say that $G$ has a $G$-marriage subject to $\mathcal{V}$. An obvious necessary condition is the {\it orbit condition}: for any $\emptyset \not = Y \subseteq [n]$, $\bigcup_{y \in Y} V_{y} \supseteq Y^{g}=\{ g(y): y \in Y \}$ for some $g \in G$. Keevash (J. Combin. Theory Ser. A 111(2005), 289--309) observed that the orbit condition is sufficient when $G$ is the symmetric group $\Sym([n])$; this is in fact equivalent to the celebrated Hall's Marriage Theorem. We prove that the orbit condition is sufficient if and only if $G$ is a direct product of symmetric groups. We extend the notion of orbit condition to that of $k$-orbit condition and prove that if $G$ is the alternating group $\Alt([n])$ or the cyclic group $C_{n}$ where $n \ge 4$, then $G$ satisfies the $(n-1)$-orbit condition subject to $\V$ if and only if $G$ has a $G$-marriage subject to $\mathcal{V}$.
\end{abstract}

\bigskip\noindent
{\sc keywords:} Hall's Marriage problem, permutation group

\section{Introduction}

In a study of the Tur{\' a}n problem for projective geometries, the following problem was first considered by Keevash, see \cite[Problem 5.1]{K}:

\noindent {\bf The $G$-Marriage Problem.}~~Let $G$ be a permutation group acting on $[n]=\{1, \ldots, n\}$ and $\V$ be a system of $n$ subsets $V_{1}, \ldots, V_{n}$ of $[n]$. When is there an element $g \in G$ so that $g(i) \in V_{i}$ for each $i \in [n]$? If such a $g$ exists, we say that $G$ has a $G$-marriage subject to $\mathcal{V}$.

An obvious necessary condition for the $G$-Marriage Problem is the {\it orbit condition} (subject to $\V$): for any $\emptyset \not
= Y \subseteq [n]$, $\bigcup_{y \in Y} V_{y} \supseteq
Y^{g}=\{ g(y): y \in Y \}$ for some $g \in G$.  Is the orbit condition also
sufficient? As noted by Keevash \cite{K}, in the case when $G$ is the symmetric group on $[n]$, the above problem is
equivalent to the Hall's marriage problem, and the necessary
and sufficient condition is that $|\bigcup_{y \in Y}V_{y}| \ge |Y|$
for every $\emptyset \not = Y \subseteq [n]$, which is equivalent to the orbit condition.

\begin{thm}\label{Hall}(Hall's Marriage Theorem)
Let $G=\Sym([n])$. Then $G$ has a $G$-marriage subject to $\V$ if and only if it satisfies the orbit condition subject to $\V$.
\end{thm}

It is natural to ask whether the orbit condition is sufficient for the $G$-Marriage Problem for other subgroups $G$ of $\Sym([n])$. One of our main results shows that the orbit condition is sufficient for the $G$-Marriage Problem if and only if $G$ is a direct product of symmetric
groups.

\begin{thm}\label{general}
Suppose $G$ is a permutation group acting on $[n]$. Then the orbit condition is sufficient for the $G$-Marriage
Problem if and only if $G$ is a direct product of symmetric groups.
\end{thm}

In view of Theorem \ref{general}, it would be interesting to find the necessary and sufficient conditions for the $G$-Marriage Problem when $G$ is not a direct product of symmetric groups. To do this, we shall require an extension of the orbit condition.

Let $k \in [n]$. We shall adopt the following notations.

\begin{itemize}
\item [(a)] $[n]^k=\overbrace{[n]\times [n]\times\cdots \times [n]}^{\textnormal {$k$ times}}$.
\item [(b)] Let $\mathbf y=(y_1,y_2,\dots, y_k)\in [n]^k$. We write $V_{\mathbf y}=V_{y_1}\times V_{y_2}\times\cdots\times V_{y_k}$.
\end{itemize}

\begin{dfn}\label{orbit-condition}
A subgroup $G$ of $\Sym ([n])$ is said to satisfy the $k$-\emph {orbit condition} subject to $\mathcal V$ if for any $\varnothing\neq Y\subseteq [n]^k$, there is a $g\in G$ such that
\begin{equation}
Y^g\subseteq\bigcup_{\mathbf y\in Y} V_{\mathbf y},\notag
\end{equation}
where $Y^g=\{ (g(y_1),g(y_2),\dots, g(y_k))\ :\ (y_1,y_2,\dots, y_k)\in Y\}$.
\end{dfn}

\noindent The following observations are obvious.

\begin{lm}\label{necessary}\
\begin{itemize}
\item [(i)] If $G$ has a $G$-marriage subject to $\mathcal V$ then $G$  satisfies the $k$-orbit condition subject to $\mathcal V$ for $k=1,2,\dots, n$.
\item [(ii)] If $G$  satisfies the $n$-orbit condition subject to $\mathcal V$ then $G$ has a $G$-marriage subject to $\mathcal V$.\hfill\qed
\end{itemize}
\end{lm}

Note that the 1-orbit condition is just the orbit condition. An example is given by Keevash \cite{K}, which shows that the 1-orbit condition is not sufficient to yield a $G$-marriage for certain group $G$. In particular, let $G$ be the subgroup generated by $(1\ 2\ 3)\in \Sym ([3])$ i.e. $G$ is a cyclic group $C_{3}$ of order 3. Clearly, $G=\{id, (1\ 2\ 3), (1\ 3\ 2)\}$, where $id$ is the identity element. Let $V_1=\{2\}$ and $V_2=\{1\}$ and $V_3=\{3\}$. It is easy to check that $G$ satisfies the 1-orbit condition subject to $\mathcal V$ but $G$ does not have a $G$-marriage subject to $\V$.

In fact, even $2$-orbit condition is not sufficient. This can be readily verified by hand or computer (we omit the details here):

\begin{prop}\label{counterexample} Let $\V$ be a system of subsets of $[3]$ consisting of $V_1=\{1,3\}$, $V_2=\{2,3\}$ and $V_3=\{1,2\}$. Then $G=\Alt([3])=C_{3}$ satisfies the 2-orbit condition subject to $\V$. However it does not have a $G$-marriage subject to $\V$.
\end{prop}

In contrast, we shall prove that the $(n-1)$-orbit condition is indeed sufficient for the $G$-Marriage Problem when $G$ is the alternating group $\Alt([n])$ or the cyclic group $C_{n}$, provided $n \ge 4$.

\begin{thm}\label{cycle} Let $G$ be the cyclic group generated by the cycle $(1\ 2\ \cdots\ n)$, $n \ge 4$. Then $G$ satisfies the $(n-1)$-orbit condition subject to $\mathcal V$ if and only if it has a $G$-marriage subject to $\V$.
\end{thm}

\begin{thm}\label{alternating} Let $G=\Alt([n])$, $n\geq 4$. Then $G$ satisfies the $(n-1)$-orbit condition subject to $\mathcal V$ if and only if it has a $G$-marriage subject to $\V$.
\end{thm}

\section{Proof of Theorem \ref{general}}

For a finite set $\Omega$, let $\Sym(\Omega)$ denote the symmetric group on $\Omega$. Suppose $G$ is a permutation group acting on $[n]$. A subset of $[n]$ is said to be a {\em base} for $G$ if its pointwise stabilizer in $G$ is trivial. The minimal size of a base for
$G$ is denoted by $\base(G)$. We refer the reader to \cite{C} for undefined terms in permutation group theory.

The most striking early result on base sizes is due to Bochert (for a survey on bases of permutation groups, see \cite{LS}):

\begin{prop}\label{Bochert}\textnormal {(Bochert \cite {Bo})}
If $G$ is a primitive permutation group of degree $n$ not
containing the alternating group $\Alt([n])$, then $\base(G) \le \frac{n}{2}$.
\end{prop}

\noindent Consequently, since $\base(\Sym([n]))=n-1$ and $\base(\Alt([n]))=n-2$, we have

\begin{prop}\label{C1}
If $G$ is a primitive permutation group of degree $n$ and $G \not = \Sym([n])$ then $\base(G) \le n-2$.
\end{prop}

\noindent Throughout this section, $\V=\{V_{1}, \ldots, V_{n}\}$ will denote a system of subsets of $[n]$ and for any subset $Y$ of $[n]$, we set

\[ V_{Y} = \bigcup_{y \in Y} V_{y}. \]

\noindent We first consider the case when $G$ is transitive.

\begin{thm}\label{transitive}
Suppose $G$ is a transitive permutation group acting on $[n]$ and the orbit condition is sufficient for the
$G$-marriage problem. Then $G = Sym([n])$.
\end{thm}

\begin{proof} Suppose that $G$ is imprimitive. Let $X_{1}, \ldots, X_{m}$ be a complete non-trivial block system
which is also a partition of $[n]$ into $m$ disjoint sets of equal size. We may assume that $X_{1} \supseteq
\{x, y\}$ and $X_{2} \supseteq \{z\}$ for some distinct elements $x, y, z \in [n]$. Construct a set system
$\V=\{V_{i}: 1 \le i \le n\}$ as follows:
\begin{eqnarray}
V_{x} & = & \{y, z\}, \nonumber \\
V_{y} & = & \{x, z\}, \nonumber \\
V_{z} & = & \{x, y\}, \nonumber \\
V_{i} & = & [n], ~~\textnormal{for all}~~i \not = x,y,z. \nonumber
\end{eqnarray}
Notice that $G$ satisfies the orbit condition subject to $\V$: let $Y \subseteq [n]$ such that $Y \not = \emptyset$. If $Y
\cap ([n]-\{x,y,z\}) \not = \emptyset$, then $V_{Y} \supseteq [n] \supseteq Y^{id}$, where $id$ is the
identity element of $G$. So we may assume that $Y \subseteq \{x,y,z\}$. But it is easy to see that if $|Y|>1$, then
$V_{Y} = \{x,y,z\} \supseteq Y^{id}$. If $|Y|=1$, then the orbit condition holds by the transitivity of $G$.

Therefore, by our assumption, there exist distinct elements $x_{i} \in V_{i}$, $1 \le i \le n$, such that the
permutation $g$, defined by $g(i)=x_{i}$, belongs to $G$. However, the image of $x, y, z$ under such a $g$ is either $y, z, x$ or $z, x, y$ respectively. In both cases, $g$ does not leave the partition $X_{1}, \ldots,
X_{m}$ invariant, which is a contradiction.

So we may assume that $G$ is primitive. Assume for a contradiction that $G \not = Sym([n])$. Let $B=\{b_{1},
\ldots, b_{k}\}$ be a minimal base of $G$ where $k=\base(G)$. Then $k \le n-2$ by Proposition \ref{C1}. Pick an element not in $B$ and denote it by $b_{k+1}$.

Construct a system $\V=\{V_{i}: 1 \le i \le n\}$ as follows:
\begin{eqnarray}
V_{b_{i}} & = & \{b_{i}, b_{k+1}\}~~\textnormal{for all}~~1 \le i \le k, \nonumber \\
V_{b_{k+1}} & = & [n]-\{b_{1}, \ldots, b_{k}, b_{k+1}\}, \nonumber \\
V_{i} & = & [n]-B ~~\textnormal{for all}~~i \not \in B \cup \{b_{k+1}\}. \nonumber
\end{eqnarray}
We now verify that $G$ satisfies the orbit condition subject to $\V$. Let $\emptyset \not = Y \subseteq [n]$. Since $k \le n-2$, all the sets $V_{i}$ are not empty and so the orbit condition holds for $Y$ when $|Y|=1$ by the transitivity of $G$. Let $|Y|>1$. Notice that if $Y \cap B \not = \emptyset$ or $b_{k+1} \not \in Y$, then $V_{Y} \supseteq Y^{id}$. So, we may suppose $b_{k+1} \in Y$ and $Y \subseteq [n]-B$. Since $|Y|>1$, we must have $V_{Y}=[n]-B$. Clearly, $V_{Y} \supseteq Y^{id}$. So the orbit condition holds.

By our hypothesis, there exists a permutation $g \in G$ such that $g(i) \in V_{i}$ for all $i \in [n]$. However, by the construction of $\V$, every such $g$ must fix $b_{1}$, $\ldots$, $b_{k}$. Since $B=\{b_{1}, \ldots, b_{k}\}$ is a base of $G$, we conclude that $g=id$. In particular, $b_{k+1} = g(b_{k+1}) \in V_{b_{k+1}}$, contradicting the fact that $b_{k+1} \not \in V_{b_{k+1}}$.
\end{proof}

\begin{prop}\label{direct-product}
The orbit condition is sufficient for the $G$-Marriage Problem if $G=\Sym(\Omega_{1}) \times \cdots \times \Sym(\Omega_{m})$.
\end{prop}

\begin{proof} Suppose $G$ satisfies the orbit condition subject to $\V=\{V_{1}, \ldots, V_{n}\}$. Then $\Sym(\Omega_{i})$ satisfies the orbit condition subject to $\V|_{\Omega_{i}}=\{V_{j} \cap \Omega_{i}: j \in \Omega_{i}\}$ for all $i=1, \ldots, m$. The result now follows immediately from Hall's Marriage theorem.
\end{proof}

\noindent {\bf Proof of Theorem \ref{general}} Assuming that the orbit condition is sufficient for the
$G$-marriage problem, we shall prove that $G$ is a direct product of symmetric groups. By Theorem
\ref{transitive}, we may suppose that $G$ is intransitive with orbits $\Omega_{i}$, $[n]=\bigcup_{i=1}^{m} \Omega_{i}$. Then $G$ is the subdirect product of its transitive
constituents $G_{1}, \ldots, G_{m}$ where $G_{i}$ is the transitive permutation group induced by the action of $G$ on the orbit $\Omega_{i}$.

Now, suppose that $G$ satisfies the orbit condition subject to $\V=\{V_{1}, \ldots, V_{n}\}$. Then, for each $i$, $G_{i}$ satisfies the orbit condition subject to $\V|_{\Omega_{i}}=\{V_{j} \cap \Omega_{i}: j \in \Omega_{i}\}$. By Theorem \ref{transitive}, we must have $G_{i}=\Sym(\Omega_{i})$ for all $1 \le i
\le m$. By Proposition \ref{direct-product}, we may assume that $G$ is not the direct product of
$\Sym(\Omega_{1}), \ldots, \Sym(\Omega_{m})$. Define the following set system $\V'$: choose a
permutation $h \in \Sym(\Omega_{1}) \times \cdots \times \Sym(\Omega_{m}) - G$. For every $1 \le i \le m$ and $j
\in \Omega_{i}$, define
\begin{eqnarray}
V'_{j} & = & \{h(j)\} \cup ([n]-\Omega_{i}). \nonumber
\end{eqnarray}
Observe that $G$ satisfies the orbit condition subject to $\V'$: let $\emptyset \not = Y \subseteq [n]$, (by the transitivity of $G_{i}$ on $\Omega_{i}$) we may assume that $Y \cap \Omega_{j} \not = \emptyset$ and $Y \cap \Omega_{j'} \not = \emptyset$ for some $j \not = j'$. Then it follows that $V'_{Y}=[n] \supseteq Y^{id}$.

Let $g \in G$ be a permutation such that $g(i) \in V'_{i}$ for all $i \in [n]$. Suppose $j \in \Omega_{i}$. Then $g(j) \in V'_{j}$. On the other hand, since $g \in G$, we must have $g(j) \in \Omega_{i}$. Therefore, $g(j) \in V'_{j} \cap \Omega_{i} \subseteq  \{h(j)\}$. Therefore, $g(j)=h(j)$ for all $j \in [n]$, i.e. $h=g \in G$, contradicting our choice of $h$.

\section{The $(n-1)$-orbit condition}

For the rest of this paper, we shall investigate the $(n-1)$-orbit condition and see when it is sufficient to yield a $G$-marriage.

\begin{lm}\label{property} Suppose $G$ satisfies the $(n-1)$-orbit condition subject to $\mathcal V$ and it does not have a $G$-marriage subject to $\V$. Then the following hold.
\begin{itemize}
\item [(a)] For each $i$, there is a $g_i\in G$ such that $g_i(i), g_i(j)\in V_j$ for all $j\neq i$ but $g_i(i)\notin V_i$.
\item [(b)] $\vert G\vert\geq n$.
\item [(c)] $V_i=[n] -\{g_{i}(i)\}$ for all $i=1,2,\dots, n$.
\item [(d)] $\{g_{1}(1),g_{2}(2),\dots, g_{n}(n)\}=[n]$.
\end{itemize}
\end{lm}

\begin{proof} (a) Let $Y=\{(i,\dots, i), (1,\dots, i-1,i+1,\dots, n)\}\subseteq [n]^{n-1}$. Then there exists a $t_i\in G$ such that $Y^{t_i}\subseteq \bigcup_{\mathbf y\in Y}V_{\mathbf y}$. There are two cases.

\noindent
{\bf Case 1.} Suppose $t_i(i)\notin V_i$. This implies that $(i,\ldots,i)^{t_{i}} \in V_{(1,\ldots,i-1,i+1,\ldots,n)}$ i.e. $t_i(i)\in V_j$ for all $j\neq i$. If $t_i(j)\in V_j$ for all $j\neq i$, then we can choose $g_i$ to be $t_i$, thus proving (a).

Suppose $t_i(j')\notin V_{j'}$ for some $j'\neq i$. Since $Y^{t_i}\subseteq \bigcup_{\mathbf y\in Y}V_{\mathbf y}$, we must have $(1, \ldots,i-1,i+1,\ldots,n)^{t_{i}} \in V_{(i,\ldots,i)}$ i.e. $t_i(j)\in V_i$ for all  $j\neq i$. This implies that $\vert V_i\vert$ has at least $n-1$ elements.

On the other hand, applying the $(n-1)$ orbit condition to the set $\{{\bf y}\}$ where $\mathbf y=(1,\dots, i-1,i+1,\dots, n)$, there exists a $g_i\in G$ such that $\mathbf y^{g_i}\in V_{\mathbf y}$. Note that $g_i(i)\notin V_i$ for $G$ does not have a $G$-marriage. So $\vert V_i\vert$ must have exactly $n-1$ elements, that is $V_i=\{t_i(1),\dots, t_i(i-1),t_i(i+1),\dots, t_i(n)\}$. Since $\{g_i(1),\dots, g_i(n)\}=[n]$, we deduce that $g_i(i)=t_i(i)$. Hence $g_i(i), g_i(j)\in V_j$ for all $j\neq i$ and $g_i(i)\notin V_i$, as desired.

\noindent
{\bf Case 2.} Suppose $t_i(i)\in V_i$.

If $(1,\ldots, i-1, i+1, \ldots,n)^{t_{i}} \in V_{(1,\ldots,i-1,i+1,\ldots,n)}$ then $t_i(j)\in V_j$ for all $j$ but this is impossible since $G$  does not have a $G$-marriage. So we must have $(1, \ldots,i-1,i+1,\ldots,n)^{t_{i}} \in V_{(i, \ldots, i)}$ i.e $t_i(j)\in V_i$ for all $j\neq i$. This implies that $V_i=\{t_i(1),\dots, t_i(n)\}=[n]$.

On the other hand, applying the $(n-1)$ orbit condition to the set $\{{\bf y}\}$ where $\mathbf y=(1,\dots, i-1,i+1,\dots, n)$, there exists a $g_i\in G$ such that $\mathbf y^{g_i}\in V_{\mathbf y}$. But $g_i(i)\in [n]= V_i$ and so $G$ has a $G$-marriage, a contradiction. Hence Case 2 cannot occur.

The proof of (a) is complete.

\noindent
(b) Now we show that all the $g_i$, $i=1,2,\dots, n$, obtained in (a) are distinct. Suppose $g_i=g_{i'}$ for some $i\neq i'$. Then $g_i(i)=g_{i'}(i)\in V_i$, a contrary to the fact that $g_i(i)\notin V_i$. Hence all the $g_i$ in (a) are distinct and $\vert G\vert\geq n$.

\noindent
(c) By (a), we see that for a fixed $k$, $V_{k}$ contains $g_i(i), g_i(k)$ for all $i\neq k$. First we show that $g_i(i)\neq g_{i'}(i')$ for $i\neq i'$. Suppose the contrary. Then $g_i(i)= g_{i'}(i')\in V_i$, contradicting the fact that $g_i(i)\notin V_i$.

So $\{g_1(1),\dots ,g_{k-1}(k-1), g_{k+1}(k+1),\dots, g_n(n)\}\subseteq V_{k}$ and $\vert V_{k}\vert\geq n-1$. We must have $\vert V_{k}\vert=n-1$, for otherwise $g_{k}(k)\in V_{k}$. Hence (c) holds.

\noindent
(d) In the proof of (c), we see that $V_{k}=\{g_1(1),\dots ,g_{k-1}(k-1), g_{k+1}(k+1),\dots, g_n(n)\}$ for all $k$ and $g_i(i)\neq g_{i'}(i')$ for $i\neq i'$. Therefore $\{g_1(1),\dots , g_n(n)\}=[n]$ and (d) holds.
\end{proof}

\begin{thm}\label{cycle_even} Let $G$ be the cyclic group generated by the cycle $(1\ 2\ \cdots\ n)$ and $n \ge 4$. Then $G$ satisfies the $(n-1)$-orbit condition subject to $\mathcal V$ if and only if it has a $G$-marriage.
\end{thm}

\begin{proof} Suppose $G$ satisfies the $(n-1)$-orbit condition subject to $\mathcal V$. Let $t=(1\ 2\ \cdots\ n)$. Then the elements in $G$ are $id, t, t^2,\dots, t^{n-1}$. Furthermore $t^u(j)=j+u ~(\textnormal {mod}\ n)$ for all $0\leq u\leq n-1$, $1\leq j\leq n$.

Assume for a contradiction that $G$ does not have a $G$-marriage. Then $G=\{g_1,g_2,\dots, g_n\}$, where the $g_i$ are as given in Lemma \ref{property}. Let $g_i=t^{u_i}$ for all $i=1,2,\dots, n$. Then all the $u_i$ are distinct. By Lemma \ref{property}, $t^{u_i}(j)\in V_j$ for all $j\neq i$, $V_{i} = [n]-\{t^{u_{i}}(i)\}$ and $\{t^{u_1}(1),t^{u_2}(2),\dots ,t^{u_n}(n)\}=[n]$.

For each $k = 1, \ldots, n$, denote the $(n-1)$-tuple $(k, k+1,\dots, n, 1, 2, \ldots, k-2)$ by ${\bf y}_{k}$. Throughout, we shall analyze the action of $G$ on these tuples. A generic $(n-1)$-tuple will be denoted by ${\bf y}$. We begin with the following claim and observation:

\noindent
{\bf Claim.} Let $\mathbf y_i^{t^u}\in V_{\mathbf y_k}$ for some $t^u\in G$. Then $V_{k-1}=[n]-\{t^u(i-1)\}$, $t^u(i-1)=t^{u_{k-1}}(k-1)$ and $u+i-1=u_{k-1}+k-1\ \textnormal {mod}\ n$.

\noindent
{\bf Proof of Claim.} Note that $t^u(i+j)\in V_{k+j}$ for $j=0,1,2,\dots,n-2$ (note that $i+j$ and $k+j$ are taken $\textnormal {mod}\ n$). This means $u+i+j\ (\textnormal {mod}\ n)\in V_{k+j}$ for $j=0,1,2, \ldots,n-2$.

Suppose $t^u(i-1)\in V_{k-1}$. Then $u+i-1\ (\textnormal {mod}\ n)\in V_{k-1}$. So $u+i+j\ (\textnormal {mod}\ n)\in V_{k+j}$ for all $0\leq j\leq n-1$. Let $u'=i-k+u \ (\textnormal {mod}\ n)$ with $0\leq u'\leq n-1$. Then $u'+k+j=u+i+j\ (\textnormal {mod}\ n)\in V_{k+j}$, which is equivalent by saying that $t^{u'} (k+j)\in V_{k+j}$. Note that $k+j$ is taken $\textnormal {mod}\ n$. So $t^{u'}(j)\in V_j$ for $j=1,2,\dots,n$ and $G$ has a $G$-marriage, a contradiction. Thus $t^u(i-1)\notin V_{k-1}$. Since $V_{k-1} = [n]-\{t^{u_{k-1}}(k-1)\}$, we conclude that $V_{k-1}=[n]-\{t^{u}(i-1)\}$, $t^{u}(i-1)=t^{u_{k-1}}(k-1)$ and $u+i-1=u_{k-1}+k-1\ \textnormal {mod}\ n$.\hfill\qed

\noindent {\bf Observation $(*)$.}~~Suppose ${\bf y}_{i}^{t^u} \in V_{{\bf y}_{k}}$ and ${\bf y}_{i'}^{t^u} \in V_{{\bf y}_{k'}}$ for some $t^u\in G$. If $i\neq i'$ then $k\neq k'$.

\noindent
{\bf Proof of Observation $(*)$.} By Claim,  $V_{k-1}=[n] -\{t^u(i-1)\}$ and $V_{k'-1}=[n]-\{t^u(i'-1)\}$. Since $i\neq i'$, $t^u(i-1) \not = t^u(i'-1)$. Therefore $k\neq k'$.\hfill\qed

\noindent
{\it Continuation of the Proof of Theorem 3.2.} For each $k=1,2,\dots ,n$, let
\begin{equation}
Y_k=\{{\bf y}_{k}, {\bf y}_{k+1}\},\notag
\end{equation}
where the subscript are taken $\textnormal {mod}\ n$. Let $k$ be fixed. The  $(n-1)$-orbit condition implies that there is a $t^u \in G$ such that $Y_k^{t_u}\subseteq \bigcup_{\mathbf y\in Y_k} V_{\mathbf y}$.

Suppose $u=u_{k-1}$. If ${\bf y}_{k}^{t^{u}} \in V_{{\bf y}_{k+1}}$, then by Claim, $t^{u_{k-1}}(k-1)=t^{u_{k}}(k)$, a contrary to the fact that $\{t^{u_1}(1),t^{u_2}(2),\dots ,t^{u_n}(n)\}=[n]$. If ${\bf y}_{k+1}^{t^{u}} \in V_{{\bf y}_{k+1}}$, then by Claim, $V_{k}=[n]-\{t^{u_{k-1}}(k)\}$, a contrary to the fact that $t^{u_{k-1}}(k)\in V_k$. Therefore ${\bf y}_{k}^{t^{u}}, {\bf y}_{k+1}^{t^{u}} \in V_{{\bf y}_{k}}$. But this contradicts Observation $(*)$.

Suppose $u=u_{k}$. If ${\bf y}_{k}^{t^{u}} \in V_{{\bf y}_{k}}$, then by Claim, $V_{k-1}=[n]-\{t^{u_{k}}(k-1)\}$, a contrary to the fact that $t^{u_{k}}(k-1)\in V_{k-1}$. If ${\bf y}_{k+1}^{t^{u}} \in V_{{\bf y}_{k}}$, then by Claim, $t^{u_{k}}(k)=t^{u_{k-1}}(k-1)$, a contrary to the fact that $\{t^{u_1}(1),t^{u_2}(2),\dots ,t^{u_n}(n)\}=[n]$.  Therefore ${\bf y}_{k}^{t^{u}}, {\bf y}_{k+1}^{t^{u}} \in V_{{\bf y}_{k+1}}$. But  this again contradicts Observation $(*)$.

Therefore, $u=u_i$ for some $i\neq k-1,k$. In particular, $t^{u_{i}}(k-1) \in V_{k-1}$ and $t^{u_{i}}(k) \in V_{k}$. By Claim, we deduce that ${\bf y}_{k}^{t^{u_i}}\in V_{{\bf y}_{k+1}}$ and ${\bf y}_{k+1}^{t^{u_i}}\in V_{{\bf y}_{k}}$.

Notice that ${\bf y}_{k}^{t^{u_i}}\in V_{{\bf y}_{k+1}}$ implies that $u_i+k-1=u_k+k\ \textnormal {mod}\ n$ (by Claim). Therefore $u_i=u_k+1 ~(\textnormal {mod}\ n)$. On the other hand, ${\bf y}_{k+1}^{t^{u_i}}\in V_{{\bf y}_{k}}$ implies that $u_i+k=u_{k-1}+k-1\ \textnormal {mod}\ n$ (by Claim). Therefore $u_{k-1}=u_i+1 ~(\textnormal {mod}\ n)$ and thus $u_{k-1}=u_k+2~(\textnormal {mod}\ n)$.

Recall that $k$ can take value $1,2,\dots ,n$. So we have the equations (in $\textnormal {mod}\ n$),
\begin{eqnarray}
u_n & = & u_1+2,  \nonumber \\
u_1 &= & u_2+2, \nonumber \\
 &\vdots &  \nonumber \\
 u_{n-1} & = & u_n+2. \label{e1}
\end{eqnarray}

Suppose $n=2m$ is even. Then by (1), $u_{n}=u_{1}+2=u_{2}+4=\cdots=u_{m}+2m=u_{m}~(\textnormal {mod}\ n)$, whence $u_n=u_m$, a contrary to the fact that all all the $u_i$ are distinct.

Suppose $n=2m+1$ is odd. Let $\mathbf z$ be the $(n-1)$-tuple $(m+1,m+1,\dots, m+1)$. Now, consider the set $A=\{{\bf y}_{1}, {\bf y}_{2}, {\bf z}\}$. Since $G$ satisfies the $(n-1)$-orbit condition, there exists $t^{u'} \in G$ such that $A^{t^{u'}} \subseteq \bigcup_{{\bf y} \in A} V_{\bf y}$.

\noindent
{\bf Case 1.} ${\bf y}_{1}^{t^{u'}} \in V_{{\bf y}_{1}}$.

By Claim, this implies that $u'+n=u_n+n\ \textnormal {mod}\ n$. Therefore $u'=u_n$. Now by Observation $(*)$, either ${\bf y}_{2}^{t^{u'}} \in V_{{\bf y}_{2}}$ or ${\bf y}_{2}^{t^{u'}} \in V_{{\bf z}}$.

Suppose ${\bf y}_{2}^{t^{u'}} \in V_{{\bf y}_{2}}$. Then by Claim,  $u'+1=u_1+1\ \textnormal {mod}\ n$, and so $u'=u_1$. But then $u_n=u_1$,  a contrary to the fact that all the $u_i$ are distinct.

Suppose ${\bf y}_{2}^{t^{u'}} \in V_{{\bf z}}$. Then $t^{u'}(j)\in V_{m+1}$ for $j=2,3,\dots, n$. Since $V_{m+1}=[n]-\{t^{u_{m+1}}(m+1)\}$, we conclude that $t^{u'}(1)\notin V_{m+1}$ and $t^{u'}(1)=t^{u_{m+1}}(m+1)$. Therefore $u'+1=u_{m+1}+m+1\  \textnormal {mod}\ n$ i.e. $u'=u_{m+1}+m\  \textnormal {mod}\ n$, and so $u_n=u_{m+1}+m\  \textnormal {mod}\ n$. On the other hand, by (1), $u_n=u_1+2=u_2+4=\cdots =u_{m+1}+2(m+1)~(\textnormal {mod}\ n)$. Therefore $u_{m+1}+m=u_{m+1}+2(m+1)=u_{m+1}+1\  \textnormal {mod}\ n$ (for $n=2m+1$). So $m-1=0\  \textnormal {mod}\ n$ whence $m=1$ and $n=3$,  a contrary to the fact that $n\geq 5$.

\noindent
{\bf Case 2.} ${\bf y}_{1}^{t^{u'}} \in V_{{\bf y}_{2}}$.

By Claim, this implies that  $u'+n=u_1+1\ \textnormal {mod}\ n$. Note that $u_1=u'-1=u'+2m~(\textnormal {mod}\ n)$. On the other hand, by (1), $u_1=u_2+2=u_3+4=\cdots =u_{m+1}+2m~(\textnormal {mod}\ n)$. Therefore $u'=u_{m+1}$.

Since $A^{t^{u'}} \subseteq \bigcup_{{\bf y} \in A} V_{\bf y}$, either $\mathbf z^{t^{u'}}\in V_{{\bf y}_{1}}$ or $\mathbf z^{t^{u'}}\in V_{{\bf y}_{2}}$ or $\mathbf z^{t^{u'}}\in V_{{\bf z}}$. If $\mathbf z^{t^{u'}}\in V_{{\bf y}_{1}}$, then $t^{u'}(m+1)\in V_{j}$ for $j=1,2,\dots ,n-1$. In particular, $t^{u'}(m+1)\in V_{m+1}$. Similarly if
$\mathbf z^{t^{u'}}\in V_{{\bf y}_{2}}$ or $\mathbf z^{t^{u'}}\in V_{{\bf z}}$, then $t^{u'}(m+1)\in V_{m+1}$. Therefore $u_{m+1}+m+1=u'+m+1~(\textnormal {mod}\ n)\in V_{m+1}$, a contradiction, for $V_{m+1}=[n]-\{t^{u_{m+1}}(m+1)=u_{m+1}+m+1~(\textnormal {mod}\ n)\}$.

\noindent
{\bf Case 3.} ${\bf y}_{1}^{t^{u'}} \in V_{{\bf z}}$.

Then $t^{u'}(j)\in V_{m+1}$ for $j=1,2,3,\dots, n-1$. Since $V_{m+1}=[n]-\{t^{u_{m+1}}(m+1)\}$, we conclude that $t^{u'}(n)\notin V_{m+1}$ and $t^{u'}(n)=t^{u_{m+1}}(m+1)$. Therefore $u'+n=u_{m+1}+m+1\  \textnormal {mod}\ n$.

Now either ${\bf y}_{2}^{t^{u'}}\in V_{{\bf y}_{1}}$ or ${\bf y}_{2}^{t^{u'}}\in V_{{\bf y}_{2}}$ or ${\bf y}_{2}^{t^{u'}}\in V_{{\bf z}}$.

Suppose ${\bf y}_{2}^{t^{u'}}\in V_{{\bf y}_{1}}$. Then by Claim, $u'+1=u_n+n\  \textnormal {mod}\ n$. Therefore $u_n=u_{m+1}+m+2\  \textnormal {mod}\ n$. On the other hand, by (1), $u_n=u_{m+1}+2(m+1)~(\textnormal {mod}\ n)$. So $u_{m+1}+m+2=u_{m+1}+2(m+1)\  \textnormal {mod}\ n$ i.e. $m=0\  \textnormal {mod}\ n$, a contradiction, for $n=2m+1$.

Suppose ${\bf y}_{2}^{t^{u'}}\in V_{{\bf y}_{2}}$. Then by Claim,  $u'+1=u_1+1\  \textnormal {mod}\ n$. Therefore $u'=u_1$ and $u_1=u_{m+1}+m+1\  \textnormal {mod}\ n$. On the other hand, by (1), $u_1=u_{m+1}+2m~(\textnormal {mod}\ n)$. So $u_{m+1}+m+1=u_{m+1}+2m\  \textnormal {mod}\ n$ i.e. $m-1=0\  \textnormal {mod}\ n$, which implies $m=1$ (for $n=2m+1$) whence $n=3$, a contradiction.

Suppose ${\bf y}_{2}^{t^{u'}}\in V_{{\bf z}}$. Then $t^{u'}(s)\in V_{m+1}$ for $s=2,3,\dots, n$. In particular, $t^{u'}(n)\in V_{m+1}$ a contrary to the fact that $t^{u'}(n)\notin V_{m+1}$.

Hence, we have shown that if $G$ satisfies the $(n-1)$-orbit condition subject to $\mathcal V$, then $G$ must have a $G$-marriage subject to $\V$.

The converse follows from Lemma \ref{necessary}.
\end{proof}

\begin{thm}\label{alternating_2} Let $G=\Alt([n])$, $n\geq 4$. Then $G$ satisfies the $(n-1)$-orbit condition subject to $\mathcal V$ if and only if it has a $G$-marriage.
\end{thm}

\begin{proof} Suppose $G$ satisfies the $(n-1)$-orbit condition subject to $\mathcal V$ and it does not have a $G$-marriage. By Lemma \ref{property},
\begin{align}
V_1& =[n]-\{a_1\}\notag\\
V_2& =[n]-\{a_2\}\notag\\
&\vdots\notag\\
V_n& =[n]-\{a_n\},\notag
\end{align}
where $\{a_1,a_2,\dots, a_n\}=[n]$. If $a_n=n$, set $p_n=id$. If $a_n\neq n$, let $b,c\in [n]- \{a_n,n\}$ and set $p_n=(a_n\ n)(b\ c)\in G$. Let $V_i^{p_n}=\{p_{n}(v) :\ v\in V_i\}$ for all $i$. Then $G$ satisfies the $(n-1)$-orbit condition subject to $\mathcal V^{p_n}$ that consist of $V_1^{p_n}, V_2^{p_n},\dots ,V_n^{p_n}$. Furthermore
\begin{align}
V_1^{p_n}& =[n]-\{a_1'\}\notag\\
V_2^{p_n}& =[n]-\{a_2'\}\notag\\
&\vdots\notag\\
V_{n-1}^{p_n}& =[n]-\{a_{n-1}'\}\notag\\
V_n^{p_n}& =[n]-\{n\},\notag
\end{align}
where $\{a_1',a_2',\dots, a_{n-1}'\}=[n-1]$.

Suppose $n-1\geq 4$. Then we can find a suitable $p_{n-1}$ as before such that $G$ satisfies the $(n-1)$-orbit condition subject to $\mathcal V^{p_np_{n-1}}$ that consist of $V_1^{p_np_{n-1}}, V_2^{p_np_{n-1}},\dots ,V_n^{p_np_{n-1}}$, and
\begin{align}
V_1^{p_np_{n-1}}& =[n]-\{a_1''\}\notag\\
V_2^{p_np_{n-1}}& =[n]-\{a_2''\}\notag\\
&\vdots\notag\\
V_{n-2}^{p_np_{n-1}}& =[n]-\{a_{n-2}''\}\notag\\
V_{n-1}^{p_np_{n-1}}& =[n]-\{n-1\}\notag\\
V_n^{p_np_{n-1}}& =[n]-\{n\},\notag
\end{align}
where $\{a_1'',a_2'',\dots, a_{n-2}''\}=[n-2]$.

Let $p=p_np_{n-1}\dots p_4$. Then $G$ satisfies the $(n-1)$-orbit condition subject to $\mathcal V^p$  with
\begin{align}
V_1^p& =[n]-\{a_1'''\}\notag\\
V_2^p& =[n]-\{a_2'''\}\notag\\
V_3^p& =[n]-\{a_3'''\}\notag\\
V_4^p& =[n]-\{4\}\notag\\
&\vdots\notag\\
V_n^p& =[n]-\{n\},\notag
\end{align}
where $\{a_1''',a_2''', a_3'''\}=[3]=\{1,2,3\}$. If $a_3=3$, set $q=id$; otherwise let $u\in [3]- \{a_3''',3\}$ and set $q=(3\ u\ a_3''')$. Then $G$ satisfies the $(n-1)$-orbit condition subject to $\mathcal V^{pq}$  with
\begin{align}
V_1^{pq}& =[n]-\{b_1\}\notag\\
V_2^{pq}& =[n]-\{b_2\}\notag\\
V_3^{pq}& =[n]-\{3\}\notag\\
V_4^{pq}& =[n]-\{4\}\notag\\
&\vdots\notag\\
V_n^{pq}& =[n]-\{n\},\notag
\end{align}
where $\{b_1,b_2\}=[2]=\{1,2\}$.

Suppose $b_1=1$. Then $b_2=2$. If $n$ is odd, then $g=(1\ 2\ \cdots \ n)\in G$ and $g(i)\in V_i^{pq}$ for all $i=1,2,\dots, n$. This implies that $p^{-1}q^{-1}g(i)\in V_i$ for all $i$ and $G$ has a $G$-marriage, a contradiction.
Suppose $n$ is even. If $n=4m$, then $g=(1\ 2)(3\ 4)\dots (4m-1\ 4m)\in G$ and $g(i)\in V_i^{pq}$ for all $i=1,2,\dots, n$. If $n=2m$ and $m$ is odd, then $g=(1\ 2\ \cdots\ m)(m+1\ m+2\ \cdots\ 2m)\in G$ and $g(i)\in V_i^{pq}$ for all $i=1,2,\dots, n$. As before $p^{-1}q^{-1}g(i)\in V_i$ for all $i$ and $G$ has a $G$-marriage, a contradiction.

Suppose $b_1=2$. Then $b_2=1$. If $n$ is odd, then $g=(1\ 3)(2\ 4)(5\ \cdots \ n)\in G$ and $g(i)\in V_i^{pq}$ for all $i=1,2,\dots, n$. But then $p^{-1}q^{-1}g(i)\in V_i$ for all $i$ and $G$ has a $G$-marriage, a contradiction.
Suppose $n$ is even. If $n=4m$, then $g=(1\ 3)(2\ 4)(5\ 6)(7\ 8)\dots (4m-1\ 4m)\in G$ and $g(i)\in V_i^{pq}$ for all $i=1,2,\dots, n$. If $n=2m$ and $m$ is odd, then $g=(1\ 3\ 5\ \cdots\ 2m-1)(2\ 4\ 6\ \cdots\ 2m)\in G$ and $g(i)\in V_i^{pq}$ for all $i=1,2,\dots, n$. But then $p^{-1}q^{-1}g(i)\in V_i$ for all $i$ and $G$ has a $G$-marriage, a contradiction.

Hence $G$ must have a $G$-marriage.

The converse follows from Lemma \ref{necessary}.
\end{proof}


\begin{thebibliography}{99}



\bibitem{Bo} A. Bochert, {\" U}ber die Zahl verschiedener Werte, die eine Funktion gegebener Buchstaben durch Vertauschung derselben erlangen kann, {\em Math. Ann.} {\bf 33} (1889), 584--590.

\bibitem{C} P. J. Cameron, Permutation groups, London Mathematical Society Student Texts {\bf 45} (1999).

\bibitem{K} P. Keevash, The T{\' u}ran problem for projective geometries, {\em Journal of Combinatorial Theory Series A}
{\bf 111} (2005), 289--309.

\bibitem{LS} M. W. Liebeck and A. Shalev, Bases of primitive permutation groups, in {\em Groups, Combinatorics and Geometry: Durham, 2001}, World Scientific, 2003.

\end{thebibliography}
\end{document}